\input amstex
\documentstyle{amsppt}
\nologo
 \pageno=1
  \loadbold
   \leftheadtext{{\smc Hongbin Chen and Yi Li}} 
    \rightheadtext{{\smc Stability for equation of Duffing type}}

\hsize=5 true in
 \vsize=8.2 true in
  \hoffset=.75 in
   \voffset= 30 pt
    \TagsOnRight
     \NoBlackBoxes

      \font\sevenrm=cmr7


\topinsert\vskip-1.5
 \baselineskip
   {\vbox{\sevenrm\baselineskip 7pt
 \noindent Manuscript submitted to\hfill Website\ http://aimSciences.org \break
  \noindent AIMS journals\hfil\break
   \line {Volume {\sevenbf X}, Number {\sevenbf X},
     XX {\sevenbf 200X}\hfil
     \eightpoint pp. X--XX}}
 }\endinsert

\footnote""{2000 Mathematics Subject Classification.
    34C25, 34C10.}
 \footnote""{Key words and phrases. Periodic solution, topological
degree, stability.}

\bigskip
 \document
  \vglue 1\baselineskip

\centerline{\bf  Existence, uniqueness, and stability of periodic
solutions}
\centerline{\bf  of an equation of Duffing type}

\bigskip
 \medskip

\centerline{\smc Hongbin Chen}
\medskip

{
 \eightpoint
  \centerline{Department of Mathematics}
   \centerline{Xi'an Jiaotong University}
    \centerline{Xi'an P.R.~China}
 } 
\medskip

\centerline{\smc Yi Li}
\medskip

{
 \eightpoint
\centerline{ Department of Mathematics }
 \centerline{ Hunan Normal University } \centerline{ Changsha, Hunan, 
 China
}
\centerline{\smc and\vrule height9pt depth5pt width0pt
}
\centerline{ Department of Mathematics }
 \centerline{ University of Iowa } \centerline{ Iowa City,
     IA 52242, 
 USA }
 } 
\medskip

\centerline{(Communicated by \qquad \qquad \qquad )}
 \medskip

{
 \eightpoint
  {\narrower\smallskip
   \noindent
    {\bf Abstract.}
 We consider a second-order equation of Duffing type.
Bounds for the derivative of the restoring force are given which
ensure the existence and uniqueness of a periodic solution.
Furthermore, the unique periodic solution is asymptotically stable
with sharp rate of exponential decay. In particular, for a restoring
term independent of the variable $t$, a necessary and sufficient
condition is obtained which guarantees the existence and
uniqueness of a periodic solution that is stable.
 \smallskip
  }
   }
    \bigskip

\def\R{{\bold R}}
\def\Re{\operatorname{Re}}
\def\diam{\operatorname{diam}}
\def\uchangedtox{x}
\def\uchangedto{\relax}
\def\refAlonso{1}
\def\refAlonsoo{2}
\def\refAmbro{3}
\def\refArnold{4}
\def\refBebernes{5}
\def\refBelitskii{6}
\def\refFabry{7}
\def\refFederer{8}
\def\refGossez{9}
\def\refGupta{10}
\def\refChen{11}
\def\refIrwin{12}
\def\refKatriel{13}
\def\refLazer{14}
\def\refLazerm{15}
\def\refLiny{16}
\def\refMawhin{17}
\def\refOrtega{18}
\def\refSell{19}
\def\refTarat{20}
\def\refStrien{21}
\def\refWolansky{22}
\def\refZitan{23}

\subheading{\S1. Introduction}
This paper is devoted to the existence, uniqueness and stability
of periodic solutions of the Duffing-type equation
$$x''+cx'+g(t,x)=h(t), \tag 1.1 $$
where $g(t,x)$ is a $T$-periodic function in $t$ and $h(t)$ is a
$T$-periodic function. The existence and multiplicity of periodic
solutions of (1.1) or more general types of nonlinear second-order
differential equations have been investigated extensively by many
authors since C.~Fabry, J.~Mawhin and M.N. Nkashama initiated the
study of the Ambrosetti--Prodi problem with periodic boundary
condition~[\refFabry]. However, the stability of periodic
solutions is less extensively studied. In [\refOrtega] R.~Ortega
studied (1.1) from a stability point of view and obtained an
Ambrosetti--Prodi-type theorem under an assumption of convex
nonlinearity. A.C. Lazer and P.J. McKenna established stability
results by converting the equation (1.1) to a fixed-point
problem~[\refLazer]. Recently, more complete results concerning
the stability of periodic solutions
of (1.1) were obtained by J.M. Alonso and
R.~Ortega~[\refAlonso,\refAlonsoo]. Under the condition that the derivative
of the restoring force is independent of $t$ and positive, they found
sharp bounds that guarantee global asymptotic stability.
In~[\refAlonso], optimal bounds for stability are obtained. But the
above results do not cover our Theorem 1, since in the theorem
the derivative of the restoring force may be negative for some $t$. From
their results, the key idea is to impose a condition on $g(x,t)$
that can rule out the existence of additional periodic solutions
that are subharmonic of order $2$. The aim of this paper is to
give conditions for existence, uniqueness, and local asymptotic
stability. The novelty of our result is not that  the friction
constant and bounds on the derivative of the restoring force
control the stability of the periodic solutions, but that the
friction constant alone determines the rate of decay  of the other
solutions of (1.1) which are near to the unique periodic solution.
More precisely, we will show that every solution to (1.1) that is
near the unique periodic solution decays uniformly at the same
exponential rate of $c/2$. Our method is based on linearization
combined with  the Floquet theory. The main results are the
following.

\proclaim {Theorem 1} Assume that $g(x,t)\in C^1({
\R}\times{\R})$, and that $g$ is $T$-periodic in $t$, such that \medpagebreak

 $(1)$ $\displaystyle\sup_{x\in \R} g'(t,x)=\alpha(t)\ll\frac{\pi^{2}}{T^{2}}+\frac{c^2}{4},$ \medpagebreak

$(2)$ there is a $T$-periodic function $\beta(t)\in C_{T}$ such
that
 $\int_0^{T}\beta(t)\,dt>0$, and $g'(t,x)\gg\beta(t)$ for all $x\in
 \R$.

Then the differential equation $(1.1)$ has a unique $T$-periodic
solution,
  which is asymptotically
  stable.

\endproclaim

 Here, we say that the periodic solution $\uchangedtox _0$ of
(1.1) is locally asymptotically stable if there exist constants
$C>0$ and $\alpha>0$ such that if $\uchangedtox $ is another solution with
$\|\uchangedtox (0)-\uchangedtox _0(0)
\|
+\|\uchangedtox '(0)-\uchangedtox _0'(0)\|=d$ sufficiently small, then
$\|\uchangedtox (t)-\uchangedtox _0(t)
\|
+\|\uchangedtox '(t)-\uchangedtox _0'(t)\|<Cde^{-\alpha t}$.
The coefficient $\alpha$ in the exponent of this upper bound is called the rate of decay of $\uchangedtox _0$.

\proclaim{Theorem 2}Assume that $g(t,x)\in C^1({
\R}\times{\R})$, such that for all $x\in \R$,
 the derivative of $g$ with respect to $x$ is subject to the bounds
$$\frac{n^2\pi^2}{T^2}+\frac{c^2}{4}\leq g'_{x}(t,x)\leq\frac{(n+1)^2\pi^2}{T^2}+\frac{c^2}{4}$$
for some  $n\geq1$.
 Then equation
$(1.1)$ has a unique $T $-periodic solution, which is stable with the
rate of decay  $c/2$.\endproclaim
 \proclaim {Theorem 3} Assume that
$g(t,x)=g(x)\in C({\R}).$ \medpagebreak

 If$$\frac{n^2\pi^2}{T^2}+\frac{c^2}{4}\leq\frac{g(x)-g(y)}{x-y}\leq\frac{(n+1)^2\pi^2}{T^2}+\frac{c^2}{4}$$
 then equation $(1.1)$ has a unique $T$-periodic solution.
Moreover, if $h \in C_{T}^{1} $ such
that
 the set
of
 critical points
 $C=\{\,t \in [\,0,T\,]:h'(t)=0\,\}$ of $h$ is Lebesgue-null,
then the periodic solution
  is
 stable with rate of decay of $c/2$.
\endproclaim
\demo{Remark} The bounds in Theorem 2 are optimal for the rate of
decay $c/2$. In the end of this section, we shall give an example
to show that the Floquet multipliers associated with (1.1) may be
a pair distinct real numbers, as soon as the derivative of the
restoring force goes a little bit across the bounds given in
Theorem 2.
\enddemo
 The following notations will be used throughout the rest of the paper.

1. $L_{T }^{p}$: $T $-periodic function $\uchangedtox \in L^{p}[\,0,T \,]$ with $%
\| \uchangedtox \| _{p}$ for $1\leq p\leq \infty $;

2. $C_{T }^{k}$: $T $-periodic function $\uchangedtox \in C^{k}[\,0,T \,]$, $k\geq
0$, with $C^{k}$-norm;

3. $\alpha(t)\gg\beta(t)$: if $
\alpha(t)\geq\beta(t)$ on $[\,0,T \,]$
and $ \alpha(t)>\beta(t)$ on some subset of \hbox{positive measure.}\!\!\!
\medskip

 We will finish this section \S1 by showing that the
bounds given in Theorem 2 are optimal for rate of decay $ c/2$.

Consider the linear equation
$$x''+cx'+q(t)x=0,
\tag 1.2$$ where $$q(t)=\cases
\frac{c^{2}}{4}+(w+\epsilon)^{2}&\text{ for } 0\leq t<\pi,
\\\frac{c^{2}}{4}+(w-\epsilon)^{2}&\text{ for } \pi\leq
t\leq\pi.\endcases$$ By the transformation $y(t)=e^{-ct/2}x(t)$,
the damping term can be eliminated and equation (1.2) reduced to
the more familiar form of the Hill equation,
$$y''+q(t)y=0.\tag 1.3$$
Evidently, the nontrivial solution $x(t)$ of (1.2) has the
rate of decay $ c/2$ if and only if $y(t)$ is a
bounded nontrivial solution of (1.3).

The fundamental solutions of the equation $y''+w^{2}y=0$ are
$x_{1}=\cos wt$, $x_{2}=\sin wt/w$. Let $T=2\pi$, the monodromy
matrix associated with (1.3) is $ A=A_{2}\cdot A_{1}$, where
$w_{1}=(w+\epsilon)^{2}$ and $w_{2}=(w-\epsilon)^{2}$, and $A_{i}$
is defined by
$$A_{i}=\pmatrix \cos\pi w_{i}&\sin\pi  w_{i}/w_{i}\\-w_{i}\sin\pi w_{i}&\cos\pi w_{i}
\endpmatrix$$
and the discriminant function is given by
$$\operatorname{tr}A=2\cos\pi w_{1}\cdot \cos\pi w_{2}
-\left(\frac{w_{1}}{w_{2}}+\frac{w_{2}}{w_{1}}\right)\sin\pi
w_{1}\sin\pi  w_{2}.$$
 The boundary of the zone of stability of (1.3) is determined by
$|\operatorname{tr}A(\epsilon)|=2$.
 According to [\refArnold] on page 120 (or by a straightforward computation)
 it can be expressed asymptotically by
 $$w=k\pm\frac{\epsilon^{2}}{k^{2}}+ o(\epsilon^{2})\text{\quad or\quad}
 w=k+\frac{1}{2}\pm\frac{\epsilon}{\pi(k+\frac{1}{2})}+o(\epsilon).$$

 If we choose $(w,\epsilon)$ in the parametric resonance region, see Figure (101) in
 [\refArnold],
   such that $(w,\epsilon)$ is near $(n/2,0)$, where $n$ is an integer,
   then one of  the Floquet multipliers of (1.3)
 is greater than one and another is less than one. In this case,
 the trivial solution of (1.2) is still asymptotically stable,
 but the rate of decay is different from $ c/2$. Thus if the
 derivative of restoring force crosses a little bit over the bounds
 given in Theorem 2, the conclusion of Theorem 2 no longer holds.
 Hence the bounds are optimal for rate of decay $ c/2$.

This paper is organized as follows. In \S2 we recall some basic
results about topological methods and prove a few lemmas that are
crucial for the proofs of the main results.  The proofs of
Theorems 1--3 are given in \S3.

\subheading{\S2. Linear periodic problems}
In this section we shall recall some basic results about
topological methods. Consider the periodic boundary value problem
$$x^{\prime }=f(t,x), \tag 2.1 $$
$$x(0)=x(T ) ,  \tag 2.2 $$
where $f\colon [\,0,T \,]\times \R^{n}\rightarrow \R^{n}$ is a continuous
function and $ T $-periodic in $t$. In order to use a homotopic method
to compute the degree, we assume that $h\colon [\,0,T\,]\times \R^{n}\times
[\,0,1\,]\rightarrow \R^{n}$ is a continuous function such that
$$h(t,x,1)=f(t,x),$$
$$h(t,x,0)=g(x),$$
where $g(x)$ is continuous. The following continuation theorem is
due to J.~Mawhin [\refMawhin].

\medpagebreak

\proclaim {Lemma 2.1}  Let $\Omega \subset C_{T }$ be an open
bounded set such that the following conditions are satisfied.

$(1)$ There is no $x\in \partial \Omega $ such that
$$x^{\prime }=h(t,x,\lambda )\;\;\;\;\;\forall\, \lambda \in
[\,0,1).$$

$(2)$  $\deg (g,\Omega \cap \R^{n},0)\neq 0.$

\noindent Then $(2.1)$--$(2.2)$ has at least one solution. \endproclaim

Let us consider the Li\'{e}nard equation
  $$x''+f(x)x'+g(t,x)=h(t), \tag 2.3$$
where $h(t)\in C_{T }$. Evidently the periodic solution of (2.3)
is equivalent to the planar system
$$x'=y-F(x),$$
$$y'=h(t)-g(t,x), \tag 2.4 $$
where $F(x)$ is a primitive of $f(x)$. A natural choice for
the homotopy in applying Lemma 2.1 is to take
$$h(t,x,y,\lambda )=\{y-F(x),\lambda h(t)+(1-\lambda
)\overline{h}-[\lambda g(t,x)+(1-\lambda )\overline{g}(x)]\},$$
where $\overline{g}(x)=\frac{1}{T}\int_{0}^{T }g(t,x)\,dt$ is
the average of $g(t,x)$. Since
$h(t,x,y,0)=(y-F(x),\overline{h}-\overline{g}(x))=G(x)$,
then the condition (2) in Lemma 2.1 reduces to
$$\deg (G,\Omega \cap \R^{2},0)\neq 0.  \tag 2.5$$
Next we consider the system (2.1) for $n=2$. We denote by $x(t,x_{0})$
the initial-value solution of (2.1) and introduce the Poincar\'{e}
map $P\colon x_{0}\rightarrow x(T ,x_{0})$. It is well known that
$x(t,x_{0})$ is a $T$-periodic solution of system (2.1) if and only if
$x_{0}$ is a fixed point of $P$. If $x$ is an isolated $T$-periodic
solution of (2.1), then $x_{0}$ is an isolated fixed point of $P$.

{\bf Definition.} A $T$-periodic solution $x$ of (2.1) will be called
a nondegenerate $T$-periodic solution if the linearized equation
$$y^{\prime }=f_{x}(t,x)y \tag 2.6$$ does not admit a nontrivial
$T$-periodic solution.

Let $M(t)$ be the fundamental matrix of (2.6) and $\mu _{1} $ and
$\mu _{2}$ be the eigenvalues of the matrix $M(T)$. Then $x(t,x_{0})$ is
asymptotically stable if and only if $\vert \mu _{i}\vert
<1$ ($i=1,2$). Otherwise, if there is an eigenvalue of $M(T )$ with
modulus greater than one, then $x(t,x_{0})$ is
unstable.

 Consider the homogeneous periodic equation
$$L_{\alpha }x=x''+cx'+\alpha (t)x=0,
\tag 2.7$$ where $c$ is constant and $\alpha (t)\in L_{T }$. From
now on, without further mention, we always suppose that the
frictional constant $c>0$.

\medpagebreak The following lemma is due to R.~Ortega [\refOrtega].\proclaim
{Lemma 2.2}Assume that $$\alpha(t)\leq
\frac{\pi^{2}}{T^{2}}+\frac{c^{2}}{4}. \tag 2.8
$$Then   equation $(2.7)$ does not admit negative Floquet multipliers.\endproclaim
 The next lemmas are crucial to the argument for
stability.
 \proclaim {Lemma 2.3} Assume
that $\alpha (t)\in L_{T }$ satisfying $(2.8)$ such that
$\overline{\alpha(t)}>0$.

  Then  the moduli of the Floquet multipliers associated with equation $(2.7)$
  are less than one.
  In other words, the trivial $T$-periodic solution is asymptotically
  stable.\endproclaim

\demo{Proof} We consider the following two cases.
\medpagebreak

Case 1. If the multipliers are a pair of conjugate numbers, the
conclusion of Lemma 2.3 follows immediately from
the Jacobi--Liouville formula.
\medpagebreak

Case 2. If the Floquet
multipliers are real numbers, then Lemma 2.2 rules out
negative multipliers, so it is
sufficient to show that the moduli of the positive multipliers are
less than one. If $x(t)$ vanishes at some $t_{0}$, then it must
vanish at $T+t_{0}$. Thus $x(t)$ solves the boundary condition
problem of the following equation:
$$(e^{ct}x')'+(\alpha(t)e^{ct})x=0, \qquad x(t_{0})=x(T+t_{0})=0. \tag2.9$$
 Next, we
consider the B.V.P.
$$(e^{ct}x')'+\lambda
e^{ct}x=0, \qquad x(t_{0})=x(t_{0}+T )=0.\tag 2.10$$
It is easy to verify that the $n$-th eigenvalue of (2.10) is
 $$\lambda_{n}=\frac{\pi^{2} n^{2}}{{T}^{2}}+\frac{c^{2}}{4}$$  and
$$y(t)=e^{-\frac{1}{2}
c(t-t_{0})}\sin\frac{\pi(t-t_{0})}{T}$$
is the first
eigenfunction corresponding to the first eigenvalue
 $$\lambda _{1}=\frac{\pi^{2}}{T^{2}}+\frac{c^{2}}{4}.$$
We compare $y(t)$ with the solution $x(t)$ of the
B.V.P. (2.9).

Since $\alpha (t)\ll \lambda _{1}$, it follows from the Sturm
comparison theorem that $y(t)$ has a zero in $(t_{0},T+t_{0})$, a
contradiction. Thus $x(t) $ does not change sign. We may assume
that $x(t) >0$,   since $x'(T)=\rho x'(0)$ and $
\frac{x'(T)}{x(T)}=\frac{x'(0)}{x(0)}$.  Dividing (2.7) by
$x(t) $ and integrating by parts gives that
$$\int_0^T\frac{x'(t)^{2}}{x(t)^{2}}\,dt+c\ln{\rho}
+\int_{0}^{T}\alpha(t)\,dt=0.$$  The hypothesis of the lemma implies
that $\ln{\rho}<0$. Hence $ 0<\rho <1$.
\enddemo
\proclaim{Lemma 2.4}Assume that there is an integer $n\geq1$ such
that
$$\frac{n^2\pi^2}{T^2}+\frac{c^2}{4}\ll
\alpha(t)\ll\frac{(n+1)^2\pi^2}{T^2}+\frac{c^2}{4}.$$ Then $(2.7)$
does not admit real Floquet multipliers.\endproclaim
\demo{Proof}If the conclusion does not hold, then there is a real
Floquet multiplier $\rho$ and a nontrivial solution $x(t)$ such
that $x(t+T)=\rho x(t)$. Since $\alpha(t)\gg\lambda_{1}$, it
follows from Sturm's separation theorem that $x(t)$ has a zero
$t_{0}\in [\,0,T\,]$. Thus $ \rho x(t_{0})=x(T+t_{0})=0$, i.e., $ x(t)
$ is a solution of the B.V.P. (2.9). The assumption of the lemma
implies that $\lambda_{n}\ll\alpha(t) \ll\lambda_{n+1}$. Therefore
$x(t)\equiv0$, a contradiction.\enddemo
 \proclaim{Lemma 2.5} Under the condition of Lemma
$2.4$, the rate of decay of any nontrivial solution of $(2.9)$
is $ c/2$.
\endproclaim
\demo{Proof}Consider the corresponding system $$X'(t)=A(t)X(t),  \tag
2.11$$ where the column vector function $X(t)=(x(t),x'(t))^{T}$
and $A(t)$ is the matrix function
$$A(t)=\pmatrix 0&1\\-p(t)&-c
\endpmatrix.$$
Let $M(t)$ be a fundamental matrix solution of
(2.11). It is well known that $M(t)$    has the form
$$M(t)=P(t)e^{Bt} \tag 2.12$$
where $P(t)$ and $B$ are $2\times2$ matrices, $P(t)=P(t+T)$ and
$B$ is constant. Let $\rho_{1}=e^{T\lambda_{1}}$ and
$\rho_{2}=e^{T\lambda_{2}}$ be the Floquet multipliers, so that
$\lambda_{1}$ and $\lambda_{2}$ are the Floquet exponents
associated with $\rho_{1}$ and $\rho_{2}$. Let $x_{1}$ and
 $ x_{2}$ be the eigenvectors of the matrix $e^{TB}$. It follows from
Lemma 2.2 that $\rho_{1}$ and    $\rho_{2}$ are
a complex conjugate pair.
Thus the eigenvectors that are associated with
different eigenvalues are linearly independent. Therefore
$y_{i}=p_{i}(t)e^{\lambda_{i}t}$ (for $i=1,2$) form the fundamental
solutions of equation (2.11). On the other hand, by applying
the Jacobi--Liouville formula we have
$${\vert\rho_{1}\vert}^{2}=\rho_{1}\rho_{2}=e^{-\int_{0}^{T}c\,dt }=e^{-cT}$$
and
$$\Re\lambda_{1}=\Re\lambda_{2}=\frac{1}{2}{\Re(\lambda_{1}+\lambda_{2})}$$
$$=\frac{1}{2T}\ln(\rho_{1}\rho_{2})
=-{\frac{c}{2}}.$$ Since every solution is a linear combination of
$y_{1}(t)$ and $y_{2}(t)$, $p_{i}(t)$ is $T$-periodic, hence it is
bounded. Therefore every nonzero solution of the equation (2.11)
decays at the same exponential rate of $c/2$.\enddemo
\medpagebreak

\subheading{\S3. Proof of main results}
Now we
prove our main results. \subsubhead Proof of Theorem 2\endsubsubhead We
divide the proof into three steps.\medpagebreak

Step 1. Existence. Without loss of generality, we may assume
that $g(0,t)=0$, for otherwise we can subtract $g(0,t)$ from both sides of
equation (1.1). Consider the parametrized equation
$$x''+cx'+{\tau}g(t,x)+(1-\tau)ax={\tau}h(t) \tag 3.1$$
for some $a\in (\lambda_{n},\lambda_{n+1} )$. We claim that there
is an $R>0$ such that equation (3.1) has no solution on $\partial
B_{R}$ for all $\tau\in[\,0,1\,]$. If there is not such an $R$, let $
x_{n}$ be a sequence such that $\|x_{n}\|\rightarrow\infty$ and
$\tau_{n}\in[\,0,1\,]$, and denote by $z_{n}$ the ratio
$\frac{x_{n}}{\|x_{n}\|}$. Dividing (3.1) by $\|x_{n}\|$, then
multiplying by $\varphi(t)\in C_{T}^{2}$ and integrating by parts,
we have that
$$\int_{0}^{T}\frac{z_{n}{\varphi}''-c{\varphi}'+[{\tau}_{n}g(t,x_{n})+(1-{\tau}_{n})a]\varphi}{\|x_{n}\|}\,dt=
{\tau}_{n}\int_{0}^{T}\frac{\varphi h_{n}}{\|x_{n}\|}\,dt. \tag
3.2$$ The condition  of Theorem 2 implies that
$\{[{\tau}_{n}g(t,x_{n})+(1-{\tau}_{n})a]x_{n}/\|x_{n}\|\}$ is
bounded. It is pre-compact in the weak* topology in $L^{1}[\,0,T\,]$.
Thus there is a subsequence such that $g(t,x_{n})/x_{n}
\rightarrow \beta(t)$ and $\tau_{n}\rightarrow\tau$. Taking the
limit in the equation (3.2), one obtains that
$$\int_{0}^{T}\{z{\varphi}''-cz{\varphi}'+w(t){\varphi}z\}\,dt=0
\tag 3.3$$ where $w(t)={\tau}\beta(t)+(1-\tau)a$. Since $w(t)$
satisfies the condition of Lemma 2.4 and $z(t)$ is a $T$-periodic
solution, it follows from Lemma 2.4 that $z(t)\equiv0$, which
contradicts $\| z(t)\|=1$. Next, by applying the homotopic
invariance property, we have that$$\deg (\phi _{0},B_{R,}0)=\deg
(\phi _{1},B_{R},0)=1.$$ This completes the proof of existence.
\medpagebreak

Step 2. Uniqueness. Let $x_{1}$ and
$x_{2}$ be two distinct $T$-periodic solutions of the equation
(1.1), and $\uchangedtox =x_{1}- x_{2}$. Then $\uchangedtox $ satisfies the equation
$$\uchangedtox ''+c\uchangedtox '+p(t)\uchangedtox =0,\tag 3.4$$ where $p(t)=[g(t,x_{1}(t))-g(t,x_{1}(t))]/[x_{1}-
x_{2}]$. Since $x(t)$ is a nontrivial $T$-periodic solution of
(3.4), it will imply that $\uchangedtox $ is an eigenfunction
associated with a Floquet multiplier equal to one. Again, Lemma
2.4 rules out this possibility. Therefore $\uchangedtox \equiv0$,
a contradiction. \medpagebreak

 Step 3. Stability. Let $x(t)$ be the unique
$T$-periodic solution obtained by Step 2. Consider the linearized
equation
$$\uchangedto v ''+c\uchangedto v '+g'(t,x(t))\uchangedto v =0 .$$ The condition of Theorem 2 implies
 that $p(t)= g'(t,x(t))$, verifying the hypotheses of Lemma 2.3, which
 gives the stability results about the linearized equation.

 To show that every solution
  of the nonlinear equation (1.1) locally decays to the unique $T$-periodic
  solution with rate of decay $c/2$,
  we need the following $C^{1}$ version of the Hartman-Grobman theorem [\refIrwin].

\medpagebreak

 \proclaim{Lemma 3.1}Let $U $ be an open neighborhood of $0$, and
  $f\colon U\subset \R^{n}\rightarrow\R^{n}$ be a $C^{1}$ function such that
 $f'_{x}(0)\colon\R^{n}\rightarrow\R^{n}$ is a contraction mapping.
 Then $f\;C^{1}$ is conjugate equivalent to
 $f'_{x}(0)$.\endproclaim

\demo{Remark}The  proof of the theorem depends  on a $C^{1}$ version of the
 Hartman-Grobman theorem. More detailed results concerning smooth
 equivalence can be found in [\refBelitskii], [\refSell], [\refStrien], [\refWolansky].  In general,
a $C^{1}$ hyperbolic map is only topologically conjugate
equivalent to the linear part in a small neighborhood of the fixed
point, thus Lemma 3.1 no longer holds for $C^{1}$ maps near the
hyperbolic fixed point without assuming that $f$ is contracting.
\enddemo

 Now we will complete the proof of Theorem 2.

  Consider the planar system associated with equation (1.1),
  $$  {\left\{\aligned x'&=y-cx,\\y'&=h(t)-g(t,x).  \endaligned \right.}\tag 3.5 $$

 Let $X_{0}(t)=(x_{0}(t),y_{0}(t))$  be the unique $T$-periodic solution determined by the initial condition
  $X_{0}(0)=(x_{0},y_{0})$. Then $X_{0}$ corresponds to the unique
  fixed point of the Poincar\'{e} mapping $PX=U(T,X)$, where
  $U(t,X)$ is the initial-value solution of (3.5) with $U(0,X)=X $. Let $ M(t) $ be
  the fundamental matrix solution of the linearization
  $$X'=A(t)X  \tag3.6$$
  of  (3.5),
  where $$A(t)=\pmatrix -c&1\\-p(t)&0
\endpmatrix.$$ By the differentiability of $X(t)$ with respect to the
initial value, the Poincar\'{e} mapping can be expressed in
terms of the initial value $X$ by the following formula:
$$PX-X_{0}=M(T)(X-X_{0})+o(X-X_{0}). \tag 3.7$$
 Referring to Lemma 2.4, $M(T)$ has a pair of conjugate eigenvalues $\lambda$, $\bar{\lambda}
 $ with $|\lambda|=e^{-cT/2}$. Thus $P(X)$ is a contracting mapping.
 According to Lemma
 3.1, there is a $C^{1}$ diffeomorphism $\varphi$ which is near enough to the
 identity that $PX-X_{0}$ is conjugate equivalent to $M(T)$.
 There is an invertible  constant matrix $C$ such that
$$ C^{-1}M(T)C=\pmatrix \lambda&0\\0&\bar{\lambda}
\endpmatrix=D(\lambda),$$
and we may suppose that $$\frac{1}{2}|X-X_{0}|<|\varphi(
X)-\varphi(X_{0})|<2|X-X_{0}|\tag 3.8 $$ for $X-X_{0}$
small, since $\varphi$ is near the
 identity. Therefore, the Liapunov exponent is given by
$$\align
\mu_{x}&=\lim_{n\rightarrow\infty}\frac{1}{n T}\ln|P^{n}X-X_{0}|
\\& =\lim_{n\rightarrow\infty}\frac{1}{n T}
\ln|\varphi\circ M(T)^{n}\circ\varphi^{-1}(X)-\varphi\circ
M(T)^{n}\circ\varphi^{-1}(X_{0})|
\\&=\lim_{n\rightarrow\infty}\frac{1}{n T}\ln|D(\lambda)^{n}\circ
C^{-1}[\varphi^{-1}(X)-\varphi^{-1}(X_{0})]|=-\frac{c}{2}.\endalign$$
The third
 equality follows from (3.8) and the fact that $C$ is invertible.

\demo{Remark}The above result shows that the Liapunov exponent is
invariant under a $C^{1}$ conjugate transformation. From the proof,
the conclusion is still true if $\varphi$ is a Lipschitz mapping.
\enddemo

Hence, the rate of decay of the solution to the unique
$T$-periodic solution  is $c/2$, independently of the initial
value $X$. Theorem 1 can be proved by exactly the same method.

 \subsubhead Proof of Theorem 3
 \endsubsubhead In order to prove Theorem 3,
  let us consider the following Li\'{e}nard equation:
  $$x''+f(x)x'+g(x)=h(t). \tag 3.9$$ The following lemma seems to be well
  known; however, we have been unable to find a published reference.
  Some similar results concerning sufficient conditions for
  existence were given in [\refGossez], [\refGupta]. We give a proof here, for
  completeness.

\proclaim{Lemma 3.2}Let $ f, g \in C(\R)$ such that $\vert
f(x)\vert>k $ for some $k>0$ and $g(x)$ is increasing. Then $(3.9)$
has a $T$-periodic solution if and only if $\overline{h}\in g(\R)$.
\endproclaim
\demo{Proof} The trivial necessary condition for existence can be
obtained by integrating (1.1) and applying the mean value theorem, and the
condition turns out to be sufficient. In fact, we consider the
parametrized equation
$$x''+f(x)x'+g(x)=\lambda h(t)+(1-\lambda) \overline{h}=h_{\lambda}(t), \tag 3.10$$
where $\lambda\in [\,0,T\,]$.
\medpagebreak

Multiplying (3.10) by $x'$ and integrating, we have that
$$k\int_{0}^{T}(x')^{2}\,dt\leq \left\vert\int_{0}^{T} x' h_{\lambda}(t)\,dt\right\vert,$$
and the H\"{o}lder inequality gives that
$$\| x'\|_{2}\leq \frac{1}{k}\| h
\|_{2}.\tag 3.11$$ On the other hand, integrating (3.10) shows that
there is a $\tau \in [\,0,T\,]$ such that $g(x(\tau))=\overline{h}$.
Since $ g(x)$ is increasing, hence $x(\tau)=g^{-1}(
\overline{h})=c $ is unique, independently of $ \lambda$, and then
a $ C^{0}$ bound can be obtained from the following
formula:$$\vert x(t) \vert =\left\vert x(\tau)+\int_{\tau
}^{t}x'(s)\,ds\right\vert \leq \vert c \vert+\frac{1}{k}\| h
\|_{2}. \tag 3.12 $$ It is easy now to obtain a $ C^{2}$ bound for
equation (3.10).\medpagebreak

Let $r_{1}$ and $r_{2}$ be
sufficiently large, and set
$$\Omega=\{\,(x,y) \mid \vert x \vert\leq r_{1},
\vert y \vert \leq  r_{2}\,\}.$$
    It follows from the estimates obtained
    above that the equivalent planar system defined in Section 2,
    $$x'=h(t,x,y,\lambda ),$$ has no solutions on $\partial\Omega $ for
    $\lambda \in [\,0,1\,]$, and the computation of the degree for $r_{1}$ and
    $r_{2}$ large enough is given by
$$\deg( G,\Omega, 0)=\deg(
g(x)-\overline{h},(-r_{1},r_{2}),0)=1.  \tag 3.13 $$ By applying
Lemma 2.1, we have that the equation (3.9) has at least one
$T$-periodic solution.\enddemo

The necessary and sufficient conditions for existence
of a
 periodic
solution can be obtained by Lemma 3.2. In order to obtain
the
stability result of Theorem 3 the following lemma concerning
the
regularity of solutions of
the
 initial-value problem
of (1.1)
 is
needed. Roughly speaking, under  Lipschitz nonlinearity, the
solution of
the
  initial-value problem
of (1.1) is
 still smooth provided
that the force term $h$ does not oscillate too violently.
\medskip

\proclaim{Lemma 3.3} Let $ u(t,\xi,\eta)$ be the solution of
the
initial-value problem
 $$  {\left\{\aligned &u''+cu'+g(u)=h(t) ,\\&u(0)=\xi,\quad u'(0)=\eta.  \endaligned \right.}\tag 3.14 $$
Assume that $g$  is a Lipschitz function,
and that
$h\in C^{1}_{T}$ such
that the set
of
 critical points
 $C=\{\,t \in [\,0,T\,]:h'(t)=0\,\}$ of $h$ is Lebesgue-null. Then, for $t\in
[\,0,\bar{t}\,]$,
 the partial derivatives of $u$ and $u'$
 with respect to $ \xi,\eta$ exist and are continuous. Moreover,
if
 $$X(t)=\pmatrix \displaystyle\frac{\partial u}{\partial \xi}&\displaystyle\frac{\partial u}{\partial \eta}
 \\\displaystyle\frac{\partial u'}{\partial \xi}&\displaystyle\frac{\partial u'}{\partial
 \eta}
\endpmatrix,$$
then $$X'(t)=A(t)X(t),\quad  X(0)=\operatorname{Id},$$
where $$A(t)=\pmatrix 0&1\\-g'(u(t))&-c
\endpmatrix.$$\endproclaim
\medskip

 \demo{Proof}
The lemma says under Lipschitz
  nonlinearity the solution of
(3.14)
is
 still smooth provided that the force term $h(t)$ does not
oscillate violently.

 We were inspired by  Lazer and McKenna's result concerning
the
 regularity of solutions of (1.1)
 for the case $g(u)=au^{+}-bu^{-}$ in [\refLazerm], for
the
 general case
of
 $g(u)$
a
 Lipschitz function; the proof
is convoluted and
 requires
 delicate analysis. We shall divide our proof into three steps,
since it involves
a great deal of
 real analysis.
\medskip

Step 1.\proclaim{Claim}
Let $u(t)$ be a solution of $(3.14)$. Let
$B\subset u([\,0,\bar{t}\,])$ be
  null.
    Then under
the
 assumptions  of Lemma $3.3$,
 the pre-image $A=u^{-1}(B)$  is null in
   $[\,0,\bar{t}\,]$.\endproclaim
   The following result concerning measure theory is needed.

\proclaim{1} Let f be differentiable on $[\,a,b\,]$,
and let
 $A$
be
a
measurable subset of $[\,a,b\,]$. If $m(f(A))=0$, then $f'(x)=0$ a.e.\
$x\in A$.\endproclaim

 Setting $$E_n=\left\{\,x\in A :\frac{\vert
f(y)-f(x)\vert}{\vert y-x\vert}> \frac{1}{n} \text{, whenever }\vert
y-x\vert<\frac{1}{n}\,\right\},$$
it is easy to check that $$E=\{x\in A: \vert f'(x)\vert
>0\}=\bigcup_{n=1}^{\infty}E_{n}.$$
In order to show that $E$ is null, it is sufficient to show that
$E_{n}$ is null for any $n$. Now, for  $n$ fixed, by
the
 additive
property of
the
 measure, it is enough to show that for any small
interval $I$ with length $\vert I\vert <\frac{1}{n}$,
the measure
 $m(
I\cap E_{n})=0
 $.  Setting $B=I\cap E_{n},$  since
$f(B)\subset f(A)$, by assumption $m(f(B))=0$.  By
the
definition of
a
null set,
there is a sequence of intervals such that
$$\bigcup_{k=1}^{\infty}I_{k}\supset f(B),\quad  \sum_{k=1}^{\infty}I_{k}<\epsilon.$$

Setting $B_{k}=f^{-1}(I_{k})\cap B$, then $B_{k}\subset
E_{n}\cap I$ and $\bigcup_{k=1}^{\infty}B_{k}=B$. Noting that
$f(B_{k})\subset I_{k}$, we have that
$$\aligned
m^{\ast}(B)&\leq \sum_{k=1}^{\infty}m^{\ast}(B_{k}) \leq
\sum_{k=1}^{\infty}\diam (B_{k}) \\&\leq \sum_{k=1}^{\infty}n \diam
(f(B_{k}))\leq n\sum_{k=1}^{\infty}m(I_{k})\leq
n\epsilon.\endaligned $$

 Since $\epsilon$ is arbitrary, it follows that $m(I\cap
E_{n})=0$. Hence, $m(E)=0$: namely, for almost all $x\in A$, we
have $f'(x)=0$. The proposition can be obtained also by directly
applying  the one-dimensional area formula. According to Theorem
3.2.5 in [\refFederer] on page 244 or in [\refLiny] on page 106,
for a given mapping,
 the area
of its image and
its
 derivative
are related by the
following formula:
  \proclaim{ 2}Let $f\colon\R\rightarrow \R$ be a Lipschitz mapping,
and for any measurable set $A$,
let
 $N(\,f\mid A,y\,)$
be the
cardinal
number of $f^{-1}(y)\cap A$.
Then  $N(\,f\mid A,y\,) \in L^{1}(\R)$,
and
$$\int_{A}\vert f'(x)\vert \,dx=\int_{\R}N(\,f\mid A,y\,)\,dy. \tag
3.15$$ \endproclaim

The cardinal number
 $N(\,f\mid A,y\,)\neq0$   if
and
 only if $y\in f(A)$ ,
so (3.15) combined with this
 assumption  gives that
$$\int_{A}\vert f'(x)\vert\,
dx=\int_{f(A)}N(\,f\mid A,y\,)\,dy=0.$$

Thus,
for almost all $x\in A$, we have $f'(x)=0$.
\medskip

We are now in a position to
finish proving the claim.\medskip

 Suppose on
the
 contrary  that $m(u^{-1}(B))>0$. It
follows from above
argument
that
 $ u'(t)=0 $ a.e.\ $t\in A=u^{-1}(B)$. We may assume that $A$
consists of accumulation points only, since the isolated points of
$A$
are
countable,
hence it is null. By Rolle's theorem, any
accumulation point of zeros of $ u'(t)$ must be a zero of $
u''(t)$. Thus,
$ u''(t)=0$ a.e.\ $  t\in A$, so the Duffing's
equation (3.14) reduces to
$$g(u(t))=h(t),\quad u'(t)=0\text{\quad a.e.\ }t\in A. \tag 3.16$$
The
 above equation gives that
$h'(t)=0$ a.e.\  $t\in A$, because
$$\frac{\vert g[u(t+s)]-g[u(t)]\vert}{\vert s\vert} \leq L \frac{\vert u(t+s)-u(t)\vert}{\vert
s\vert},$$ where $L$ is the Lip constant of $g$. This contradicts
our assumption that the set of critical points of $h$ is null
because of $m( A)>0$, which completes our claim.
\medskip

Step 2. The existence of partial derivatives.

 \medpagebreak

 Let $L$ be the Lipschitz constant of $g$ and let $
\{\xi_{n}\}_{1}^{\infty}$ be a sequence
 of numbers such  that
 $\xi_{n}\rightarrow \xi_{0}$ as $n\rightarrow \infty $. Let
   $$u_{n}=u(t,\xi_{n},\eta),\qquad u_{0}=u(t,\xi_{0},\eta)$$ and let $$\Psi_{n}=(u_{n}-u_{0})/(\xi_{n}-\xi_{0}).$$
 Since $u_{n}''+cu_{n}'+g(u_{n})=h(t)$ for $n\geq0$, we have that
 $$ \vert\Psi_{n}(t)''\vert\leq c\vert\Psi_{n}(t)'\vert+L\vert\Psi_{n}(t)\vert. \tag 3.17$$
 Moreover, we have that $\Psi_{n}(0)=1$, $\Psi_{n}'(0)=0$, so for $n\geq 1$ and $t\geq0$,
 $$\Psi_{n}(t)=1+\int_{0}^{t}\Psi_{n}(t)'\,dt,\qquad\Psi_{n}'(t)=1+\int_{0}^{t}\Psi_{n}(t)''\,dt.\tag3.18$$
 Let $M=1+c+L$; then (3.17) gives that
 $$ \vert\Psi_{n}(t)\vert+ \vert\Psi_{n}'(t)\vert\leq 1+M\int_{0}^{t}\vert\Psi_{n}(s)\vert+ \vert\Psi_{n}'(s)\vert \,ds.$$
 It follows from Gronwall's inequality that
 $$ \vert\Psi_{n}(t)\vert+ \vert\Psi_{n}'(t)\vert\leq\exp(Mt).\tag 3.19$$
{}From (3.17) and (3.19) it follows that the sequences
$\{\Psi_{n}(t)\}$ and $\{\Psi_{n}'(t)\}$ are equicontinuous and
uniformly bounded on $[\,0,\bar{t}\,]$, so there exists a
subsequence [still denoted by $\{\Psi_{n}(t)\}$] and a $C^{1}$
$z(t)$ such that $\Psi_{n}(t)\rightarrow z$ and
$\Psi_{n}'(t)\rightarrow z'$ as $n\rightarrow \infty $ uniformly
on $[\,0,\bar{t}\,]$. Moreover, $\{\Psi_{n}(t)\}$ satisfies the
following equation:
$$\Psi_{n}(t)''+c\Psi_{n}(t)'+\frac{[g(u_{n}(t))-g(u_{0}(t))]}{\xi_{n}-\xi_{0}}=0. \tag3.20$$
Since $g$ is a Lipschitz function,
it follows that
 $g'(u)$ exists for almost every
$u$. Let $E$ be a set such that
the
 derivative of $g$ does not exist;
then $m(E)=0$. Since $\Psi_{n}(t)\rightarrow z$, thus, the limits
of the last term in (3.20) exist for $u_{0}(t)\in
u_{0}([\,0,\bar{t}\,])\setminus E$. Let $A=u_{0}^{-1}(E)$ be a
pre-image of $E$. It follows from step 1 that $A$ is null, namely,
for almost all $t \in [\,0,\bar{t}\,]$ the limit
$$\lim_{n\rightarrow\infty}\frac{[g(u_{n}(t))-g(u_{0}(t))]}{\xi_{n}-\xi_{0}}=g'(u_{0}(t))
z(t),$$and
the
 Lipschitz condition gives that
the
 above sequence is
bounded. Taking the limit in (3.18), it follows from the
 Lebesgue  bounded
dominated convergence theorem that$$z(t)=1+\int_{0}^{t}z(s)\,ds,\quad
z'(t)=-\int_{0}^{t}cz'(s)+g'(u_{0}(s)) z(s)\,ds.$$

Thus, $z'(t)$ is precisely
the
function satisfying
$$z''+cz'+g'(u_{0}(t))z=0,\quad z(0)=1,\quad z'(0)=0. $$
Since this determines $z$ uniquely,
the
 original sequences
$\{\Psi_{n}(t)\}$ and
$\{\Psi_{n}'(t)\}$ must converge to $z$ and $z'$ respectively,
which gives the existence of $\partial u/\partial\xi$ and
$\partial u'/\partial\xi$. Similarly, we can show that $\partial
u/\partial\eta$ and $\partial u'/\partial\eta$
exist.

\medpagebreak

Step 3. Continuity of partial derivatives.

\medpagebreak

In order to show
the
 continuity of
the
 partial derivatives with
respect to
the
 initial value, it is sufficient to show that if $
\{\xi_{n}\}_{1}^{\infty}$, $ \{\eta_{n}\}_{1}^{\infty}$
are sequences such that $\xi_{n}\rightarrow \xi_{0}$,
$\eta_{n}\rightarrow \eta_{0}$, if
$u_{n}=u(t,\xi_{n},\eta_{n})$, $u_{0}=u(t,\xi_{0},\eta_{0})$
 and if
 $$y_{n}''+cy_{n}'+g'(u_{n}(t))y_{n}=0,\quad y_{n}=1,\quad y_{n}'(0)=0,
$$
then $y_{n}\rightarrow y_{0}$, $y_{n}'\rightarrow y_{0}'$
uniformly on $[\,0,\bar{t}\,]$.

\medpagebreak

First, let us show that $g'(u_{n}(t))\rightarrow g'(u_{0}(t)) $ in
 $L^{1}$.

 Since $g'(u)=\lim_{n\rightarrow
\infty}n[g(u+\frac{1}{n})-(g(u))]$ is the limit of
a
 continuous
function,
it follows that
 $g'(u)$ is  a measurable function on
$X=u_{0}([\,0,\bar{t}\,])$. By Lusin's theorem, [\refFederer] on
page 76, for every $\delta>0$ there is a bounded closed subset
$E_{1}$ of $X=u_{0}([\,0,\bar{t}\,])$ with $m(X\setminus
E_{1})<\delta$ such that $g'(u)$ is continuous on $E_{1}$. Again,
by applying Lusin's theorem to the subset $(X\setminus E_{1})$,
there is a closed set $E_{2}$ of $(X\setminus E_{1})$ with
$m(X\setminus (E_{1}\cup E_{2}))<\delta/2$ such that $g'(u)$ is
continuous on $E_{2}$. Evidently, $E_{1}\cap E_{2}=\varnothing$.
Thus, the distance $d(E_{1}, E_{2})>0$, which implies that the two
sets
 do not have
a
 limit
point in common, so $g'(u)$ is continuous on $E_{1}\cup
E_{2}$. By repeating this process, we obtain a sequence of closed
subsets $F_{k}$ of $X$ such that $F_{1}=E_{1}$,
$F_{2}=E_{1}\cup
E_{2}$, $\dots$,
$F_{k}= \bigcup_{j=1}^{k}E_{j}$, with
the
 following
properties:
1)~$g'(u)$ is continuous on $F_{k}$; 2)~$F_{k}\subset
F_{k+1}$, $m(X \setminus F_{k})\rightarrow 0$, $k\rightarrow \infty$.
Next,
set $$G_{k}=\{\,t\in [\,0,\bar{t}\,] :u_{0}(t,\xi_{0},\eta_{0})\in
X\setminus F_{n} \,\}.$$

We claim that $m(G_{k})\rightarrow 0$ as $n\rightarrow \infty$.
If not, since $G_{k}\supset G_{k+1}$, there is an
$\varepsilon_{0}>0$ such that $m(G_{n})>\varepsilon_{0}$. By
the
monotone property of the measure, we have that
$$m\left(\bigcap_{k=1}^{\infty}G_{k}\right)=\lim_{k\rightarrow\infty}m(G_{k})>\varepsilon_{0}.$$

But the image of $A=\bigcap_{k=1}^{k=\infty}G_{k}$ under $u_{0}$
is
contained in $\bigcap_{k=1}^{k=\infty} (X\setminus F_{k})$, which
is a null subset of $X$. According to step 1, $m(A)=0$, which leads
to a
contradiction.

It follows from
the
 theorem concerning continuous dependence on
the
initial value that $u_{n}\rightarrow u_{0}$  uniformly on
$[\,0,\bar{t}\,]$. Now
since
$m(G_{k})\rightarrow 0$ as $n\rightarrow
\infty$,
we may
 choose $k$ large enough that $m(G_{k})<\epsilon$. Now we
fix $k$,
and
 since $ g'$
depends
 continuously on $F_{k}$,
it follows that
$g'(u_{n}(t))\rightarrow g'(u_{0}(t))$ uniformly on
$G_{k}^{c}=[\,0,\bar{t}\,]\setminus G_{k}$, which
in turn
 implies that
$$\aligned &\lim_{n\rightarrow \infty}\int_{0}^{\bar{t}}\vert g'(u_{n}(t))-
g'(u_{0}(t))\vert \,dt\\&\qquad\qquad=\lim_{n\rightarrow \infty}\int_{G_{k}}\vert
g'(u_{n}(t))- g'(u_{0}(t))\vert \,dt+\int_{G_{k}^{c}}\vert
g'(u_{n}(t))- g'(u_{0}(t))\vert \,dt \\&\qquad\qquad\leq 0+2L\epsilon.\endaligned
$$

Therefore, $g'(u_{n}(t))\rightarrow g'(u_{0}(t)) $ in  the sense
of $L^{1}$.

The rest of
the
 proof is similar to
that of
 Lemma 2.2 in [\refLazerm]. Setting
$v_{n}(t)=y_{n}(t)-x_{n}(t)$, then
$${\left\{\aligned &v_{n}''+cv_{n}'+g'(u_{n}(t))v_{n}=[g'(u_{0}(t))-g'(u_{n}(t))]v_{0}',
\\& v_{n}(0)=v_{n}'(0).\endaligned \right.} \tag 3.21$$
Denote the $L^{1} $ norm of $g'(u_{n}(t))- g'(u_{0}(t))$ by
$\varepsilon_{n}$ and
the
 $C^{0}$ norm of $v_{0}$ by $k$. Then from
(3.21), we have that
 $$\vert v_{n}(t)\vert+\vert v_{n}'(t)\vert \leq \varepsilon_{n}k +M
 \int_{0}^{t}[\vert v_{n}(t)\vert+\vert v_{n}'(t)\vert] \,dt,$$
 where $M=1+c+L$.
The Gronwall inequality implies  that
 $$\vert v_{n}(t)\vert+\vert v_{n}'(t)\vert \leq \varepsilon_{n}k
 e^{Mt}.$$

This completes
the proof
  that $y_{n}\rightarrow y_{0}$,
$y_{n}'\rightarrow y_{0}'$
uniformly on $[\,0,\bar{t}\,]$.

 \enddemo

On the same lines as the proof
of
 Theorem 2, the stability
result of Theorem
3 follows from Lemma 3.3.
\medskip

 \noindent {\bf Acknowledgment:}\enspace We wish to thank
the referee for his/her helpful comments and suggestions.

\medpagebreak

\medskip

E-mail address : hbchen\@mail.xjtu.edu.cn, yi-li\@uiowa.edu

\medskip

Received XX; revised XX.


\Refs \ref\no\refAlonso \by J.M. Alonso\paper Optimal intervals of
stability of a forced oscillator\jour Proc. Amer. Math. Soc. \vol
123\yr1995\pages 2031--2040\endref

\ref\no\refAlonsoo \by J.M. Alonso and R.~Ortega\paper Boundedness
and global asymptotic stability of a forced oscillator\jour
Nonlinear Anal. \vol25\yr1995\pages 297--309\endref

\ref\no\refAmbro \by A.~Ambrosetti and G.~Prodi \book A Primer of
Nonlinear Analysis \publ Cambridge University Press \publaddr
Cambridge \yr1993\endref

\ref\no\refArnold \by V.I. Arnol$'$\kern-1.5ptd  \book Mathematical Methods of Classical Mechanics
\bookinfo Graduate Texts in Mathematics, vol.~60\publ Springer-Verlag \yr1989
\endref

\ref\no\refBebernes \by J.~Bebernes and
M.~Martelli\paper Periodic solutions for Li\'{e}nard systems \jour Boll. Un. Mat. Ital. (5)
\vol 16 \yr1979\pages 398--405\endref

\ref\no\refBelitskii\by G.R. Belitskii
\paper  Functional equations, and conjugacy of local diffeomorphisms of finite smoothness class
\jour Functional Anal. Appl.\vol7\yr1973\pages 268--277\endref

\ref\no\refFabry \by C.~Fabry, J.~Mawhin, and
M.N. Nkashama \paper A multiplicity result for periodic solutions of forced nonlinear second
order ordinary differential equations\jour Bull. London Math. Soc.\vol18 \yr
1986\pages173--180\endref

\ref\no\refFederer\by H.~Federer \book Geometric Measure Theory
\bookinfo Grundlehren Math. Wiss., vol.~153
\publ Springer-Verlag \publaddr Berlin-Heidelberg \yr 1969
\endref

\ref\no\refGossez\by J.-P. Gossez and P.~Omari \paper Periodic
solutions of a second order ordinary differential equation: a
necessary and sufficient condition for nonresonance \jour
J.~Differential Equations \vol 94 \yr 1991 \pages 67--82 \endref

\ref\no\refGupta\by C.P. Gupta, J.J. Nieto, and L.~Sanchez \paper
Periodic solutions of some Li\'{e}nard and Duffing equations \jour
J.~Math. Anal. Appl. \vol 140 \yr 1989 \pages 67--82\endref

\ref\no\refChen \by H.B. Chen, Yi Li, and X.J. Hou \paper Exact
multiplicity for periodic solutions of
 Duffing type \jour Nonlinear Anal.
\vol55\yr2003\pages115--124
\endref

\ref\no\refIrwin \by M.C. Irwin \book Smooth Dynamical Systems
\bookinfo Pure and Applied Mathematics, vol.~94\publ Academic
Press\publaddr London \yr 1980\endref

\ref\no\refKatriel \by G.~Katriel \paper Uniqueness of periodic solutions
for asymptotically linear Duffing equations with strong forcing
\jour Topol. Methods Nonlinear Anal. \vol 12 \yr
1998\pages263--274
\endref

\ref\no\refLazer\by A.C. Lazer and P.J. McKenna \paper On the
existence of stable periodic solutions of differential equations
of Duffing type\jour Proc. Amer. Math. Soc.\vol110
\yr1990\pages125--133\endref

\ref\no\refLazerm\by A.C. Lazer and P.J. McKenna
\paper
Existence, uniqueness, and stability of oscillations in
differential equations with asymmetric nonlinearities
\jour Trans. Amer. Math. Soc.
\vol 315 \yr 1989 \pages 721--739
\endref

\ref\no\refLiny\by Lin Fanghua and Yang Xiaoping\book Geometric
Measure Theory---An Introduction
\bookinfo Advanced Mathematics (Beijing/Boston), vol.~1
\publ International Press
\publaddr Boston\yr 2002
\endref

\ref\no\refMawhin \by J.~Mawhin \paper  Topological degree and
boundary value problems for nonlinear differential equations
\inbook  Topological Methods for Ordinary Differential Equations
(Montecatini Terme, 1991) \eds M.~Furi and P.~Zecca \publ Lecture
Notes in Math., vol.~1537, Springer\publaddr Berlin \yr1993 \pages
74--142\endref

\ref\no\refOrtega \by R.~Ortega\paper Stability and index of
periodic solutions of an equation of Duffing type\jour Boll. Un.
Mat. Ital. B (7)\vol3 \yr 1989 \pages 133--146\endref

\ref\no\refSell\by G.R. Sell
\paper Smooth linearization near a fixed point
\jour Amer. J. Math \vol107 \yr1985\pages1035--1091\endref

\ref \no\refTarat \by G.~Tarantello\paper  On the number of solutions for the forced pendulum equation
\jour J.~Differential Equations \vol80\yr 1989\pages 79--93\endref

\ref\no\refStrien\by S. van Strien
\paper Smooth linearization of hyperbolic fixed points without
              resonance conditions
\jour J.~Differential Equations \vol85\yr1990\pages66--90\endref

\ref\no\refWolansky \by G. Wolansky
\paper Limit theorem for a dynamical system in the presence of
              resonances and homoclinic orbits
\jour J.~Differential Equations\vol83\yr1990\pages300--335\endref

\ref\no\refZitan \by A.~Zitan and R.~Ortega \paper Existence of
asymptotically stable periodic solutions of a forced equation of
Li\'{e}nard type\jour Nonlinear Anal.\vol
22\issue8 \yr1994 \pages993--1003\endref \endRefs
\enddocument